\documentclass[reqno,tbtags,intlimits,a4paper,oneside,11pt]{amsart}
\usepackage{amsfonts,amsmath,amssymb}
\newtheorem{Theorem}{Theorem}[section]
\newtheorem{prop}[Theorem]{Proposition}

\DeclareMathOperator{\arcosh}{arcosh}

\def\beq#1#2\eeq{%
        \begin{equation}%
        \label{#1}%
            #2%
        \end{equation}%
    }

\usepackage{color}
\usepackage{graphicx}
\usepackage{epstopdf}

\title[Fricke, Frobenius and Markov]{Fricke identities, Frobenius $k$-characters and Markov equation}

\author{V.M. Buchstaber}\address{Steklov Mathematical Institute and Moscow State University, Russia}
\email{buchstab@mi-ras.ru}

\author{A.P. Veselov}
\address{Department of Mathematical Sciences,
Loughborough University, Loughborough LE11 3TU, UK  and Moscow State University, Moscow 119899, Russia}
\email{A.P.Veselov@lboro.ac.uk}

\begin{document}

\maketitle

\maketitle
\vspace{0.5cm}
\rightline{\em To the memory of our dear friend Boris Dubrovin}
\bigskip

\begin{abstract}
In 1896 Frobenius and Fricke had published two seemingly unrelated papers: Frobenius had started to develop his theory of $k$-characters for finite groups motivated by Dedekind's question about factorisation of the group determinant, while Fricke followed Klein's approach to the uniformization theorem. We show that in fact these two works can be naturally linked and both are related to remarkable Markov's paper of 1880 on arithmetic of binary quadratic forms.
\end{abstract}

%\tableofcontents

\section{Introduction}

In 1880 Andrei A. Markov, St Petersburg student of Korkin and Zolotarev, in his Master Thesis  \cite{Markov} discovered a remarkable relation between the theory of binary quadratic forms and the following Diophantine equation known as {\it Markov equation}
\begin{equation}
\label{Markov}
x^2 + y^2+ z^2 = 3xyz.
\end{equation}
Markov showed that all positive solutions of this equation can be found from an obvious solution (1,1,1) by the Vieta involution 
\begin{equation}
\label{Markovinv}
\tau: (x,y,z)\rightarrow (x,y, 3xy-z)
\end{equation}
and permutations.
These celebrated {\it Markov triples} play an important role in Diophantine Analysis: as it became clear after the work of Hurwitz \cite{Hurwitz} the corresponding quadratic irrationals 
$$
\alpha=\frac{b}{x}+\frac{y}{xz}-\frac{3}{2}+\frac{\sqrt{9z^2-4}}{2z}, \quad by-ax=z,
$$
are in some precise sense the most irrational numbers (see more details in \cite{Aigner, Delone}).

In 1896 there were appeared two papers \cite{Fricke, Frob1} by Fricke and Frobenius respectively, which seemingly had nothing to do neither with the subject of Markov's work, nor with each other:
Fricke studied the uniformization problem for one-holed tori, while Frobenius started the development of the group character theory. 

However, starting from 1950s work of Gorshkov and Cohn \cite{Cohn, Gorshkov}, it was realised that Fricke's work is directly related to Markov's study \cite{Haas, Series}. Nowadays, all this is simply a part of the Teichm\"uller theory of one-holed tori \cite{Goldman, Keen1971, Keen}. 

On the contrary, a close relationship of this circle of problems with the work of Frobenius was not noticed until very recently, when we have realised this while working on our joint paper \cite{BV}. 

The aim of the present article is to explain an unexpected link between the works of Frobenius and Fricke (and thus with Markov's work) in more detail. 

The central part is a brief review and extension of the Frobenius theory of $k$-characters for finite groups, which was motivated by the factorization problem of the group determinant. We closely follow his original work \cite{Frob2}, which we found very deep and still worthy to be fully appreciated.

We should mention that recently Markov triples attracted a lot of interest because of the links with the theory of Frobenius manifolds and related Painlev\'e-VI equation \cite{D1} and algebraic geometry \cite{HP, Rud}. In fact, one of us (APV) became interested in Markov equation after learning from Boris Dubrovin about his discovery of their remarkable link with the enumerative geometry and quantum cohomology of $\mathbb P^2$.

We dedicate this paper to the memory of Boris Dubrovin, our dear friend and colleague, whom we will always miss a lot.

%From 1950s work by Gorshkov \cite{Gorshkov} and Cohn \cite{Cohn} this story became part of the hyperbolic geometry and the Teichm\"uller space theory, see e.g. \cite{Fock, Goldman, Malyshev}.
%
%The main idea goes back to the works of Fricke and Klein on uniformisation problem for the punctured tori $T_*^2$. To solve this problem we need to define a suitable homomorphism of $\pi_1(T_*^2)$ into the group of motion $PSL_2(\mathbb R)$ of the hyperbolic plane $\mathbb H^2.$ 
%
%Since the fundamental group $\pi_1(T_*^2)=F_2$ is free group with two generators, we need to find two hyperbolic matrices $A$ and $B$ from  $SL_2(\mathbb R)$ such that their commutator $ABA^{-1}B^{-1}$ is parabolic (which is the puncture condition). One can show \cite{Beardon} that this commutator must have two eigenvalues equal to -1, so we must have
%$$tr\, (ABA^{-1}B^{-1})=-2.$$
%

\section{Fricke identities}

Robert Fricke's paper \cite{Fricke} is a nice read although it is written in a style, typical for that time, without identifying statements as lemmas, theorems etc. 
We present now his findings using modern notations and terminology.

The main subject of the paper is the uniformization problem for the one-holed tori. The corresponding hyperbolic structures are the quotients $\mathbb H^2/G$, where $G \subset SL_2(\mathbb R)$ are discrete subgroups of $SL_2(\mathbb R)$, which are called after Poincare {\it Fuchsian groups.}\footnote{It is well-known that Klein did not like this terminology and criticised Poincare for introducing it, which is probably why Fricke, being a student of Klein, did not use it.}
The fundamental group of a one-holed torus is free group with two generators, so the corresponding subgroup $G$ is freely generated by two matrices $A,B \in SL_2(\mathbb R).$

Fricke starts with the following identity, which is true for any $A,B \in SL_2(\mathbb R)$ and $C=AB:$
\beq{fricke}
(tr\, A)^2+(tr\, B)^2+(tr\, C)^2=tr\, A \, tr\, B\, tr\,C +tr\, (ABA^{-1}B^{-1})+2
\eeq
(see formula (6) in \cite{Fricke}). 

Fricke did not give the proof, but this is elementary. Indeed, by Cayley-Hamilton theorem any matrix $Y\in SL_2(\mathbb R)$ satisfies its characteristic equation $\lambda+\lambda^{-1}=tr\, Y,$ which means that
$
Y+Y^{-1}=(tr\, Y) I.
$
Multiplying it by $X$ we have
$
XY+XY^{-1}=(tr\, Y) X
$
and, after taking trace, the well-known identity (sometimes called as skein relation)
\beq{fricke0}
tr\, XY+tr\, XY^{-1}=tr\, X \, tr\, Y
\eeq
valid for all $X,Y \in SL_2(\mathbb R).$ The Fricke identity (\ref{fricke}) is the result of several applications of this identity. 

Indeed, applying (\ref{fricke0})  to $X=ABA^{-1}$ and $Y=B$ we have
$$
tr\, ABA^{-1}B+tr\, ABA^{-1}B^{-1}=tr\, ABA^{-1} \, tr\, B=tr\, B^2.
$$
Applying now (\ref{fricke0}) to $X=AB$ and $Y=A^{-1}B$ we have
$$
tr\, ABA^{-1}B+tr\, ABB^{-1}A=tr\, AB \, tr\, A^{-1}B.
$$
Since $AB=C$ and $tr\, A^{-1}B=-tr\, AB+tr\, A \, tr\, B$ (again by (\ref{fricke0}), we have
$$
tr\, ABA^{-1}B=-tr\, (A^2)-(tr\, C)^2+tr\, A \, tr\, B\, tr\,C.
$$
Using now that $tr\, (A^2)=(tr\, A)^2-2$ for any $A\in SL_2(\mathbb R)$, we get the Fricke identity (\ref{fricke}).

Denoting $x:=tr\,A, \, y:=tr\,B, \,z:=tr\,C \,$ and $j_c=tr\, ABA^{-1}B^{-1},$ Fricke arrived at the real surface determined by the following cubic equation
\beq{fricke2}
x^2+y^2+z^2-xyz -(j_c+2)=0
\eeq
with positive $x,y,z$ (see formula (8) in \cite{Fricke}).

Next important observation of Fricke is that the natural action of the modular group $SL_2(\mathbb Z)$ on this surface is generated by two transformations: in Fricke's notations
\beq{ST}
S: x'=y,\, y'=x,\, z'=xy-z, \quad \quad
TS: x'=y, \,y'=z, \,z'=x.
\eeq
If we add all permutation of variables, then we have the action of $GL_2(\mathbb Z)$ (which Fricke called extended modular group).

Fricke finishes with an explicit geometric description of the group $G$ corresponding to a positive real triple $(x,y,z)$ satisfying (\ref{fricke2}).

Fricke did not make any connections to earlier Markov's results \cite{Markov}, because his topic was completely different and had nothing to do with number theory.

However, if we consider the pairs of hyperbolic matrices $A,B$ with parabolic commutator $ABA^{-1}B^{-1}$, then $j_c=tr\, ABA^{-1}B^{-1}=-2,$ then Fricke's relation (\ref{fricke2}) will take form 
\beq{fricke3}
x^2+y^2+z^2-xyz=0.
\eeq
If we rescale the variables by $1/3$ and assume that they are integer, we will come to Markov equation (\ref{Markov}).

From modern point of view the condition $j_c=-2$ corresponds to the {\it punctured tori}, so the surface given by (\ref{fricke3}) with $x,y,z>0 $ can be interpreted as Teichm\"uller space of such tori. The positive real component of the cubic surface given by Fricke relation (\ref{fricke2}) is the Teichm\"uller space of one-holed tori with geodesic boundary of length
 $l=2\arcosh (-\frac{j_c}{2}).$

This fact was probably first time explicitly mentioned by Keen \cite{Keen1971} (see also Goldman \cite{Goldman}), but the relation of Markov's results with hyperbolic geometry was discovered much earlier in the 1950s independently by Gorshkov and Cohn \cite{Cohn, Gorshkov}. 

%Gorshkov and Cohn showed that Markov's $PGL_2(\mathbb Z)$-orbit  corresponds to the punctured $\mathbb Z_3$-symmetric rhombic (equianharmonic) torus with complete hyperbolic metric: Markov numbers can be interpreted as $\frac{2}{3} \cosh l$ the lengths $l$ of the corresponding simple closed geodesics. Markov triples describe the lengths of  such geodesics with one-point pairwise intersections \cite{Cohn, Gorshkov, Haas, Series}.

More precisely, consider the punctured equianharmonic torus $\mathbb T^2_*$ with conformally equivalent complete hyperbolic metric, then Markov triples $(x,y,z)$ are simply related to the lengths of the triples of simple geodesics on this particular punctured torus $\mathbb T^2_*$ with pairwise intersection numbers 1 (see more details in Haas \cite{Haas} and Series \cite{Series}).

\section{Frobenius $k$-characters}

Ferdinand Georg Frobenius had two relevant and both very important papers \cite{Frob1, Frob2}, written in 1896. Prior that he had a letter exchange with Richard Dedekind, his predecessor at ETH Zurich. Dedekind informed him about his unpublished results and conjectures on the factorisation of the group determinant, which Frobenius immediately got keen interest in.
Dedekind's motivation came from number theory, in particular by the results of Dirichlet.

It took five months of 1896 for Frobenius to completely resolve Dedekind problem of the factorisation problem for group determinant. For this Frobenius had to develop the character theory for finite groups, which was created just for this purpose (and got its proper status after the work by Issai Schur, a brilliant student of Frobenius).

We agree with Slava Gerovitch, who wrote in nicely written historic interludes in \cite{Etingof} about these two papers of Frobenius: "Although his proofs have now be supplanted by easier modern version, his skills as group determinant virtuoso remain unsurpassed, as the tool itself went out of use."

What is much less known however, is that Frobenius developed also the theory of the so-called $k$-characters, which we are going to describe now in more details.

Let us start with explaining the Dedekind's group determinant. Let $G$ be a finite group of order $n$ and introduce the variables $x_g$ labelled by the elements $g\in G$. Let's order them as $g_1=e, g_2, \dots, g_N$ and consider the multiplication $N\times N$ matrix $X_G$ with $X_{i,j}=x_{g_ig_j^{-1}}.$

The {\it group determinant} of $G$ is the determinant of $X_G$:
\beq{gdet}
\Theta_G=\det \, X_G.
\eeq
For the cyclic group $G=\mathbb Z_N$  the group determinant is well-known {\it circulant.} For example, for $G=\mathbb Z_3$
$$
X=\begin{pmatrix}
  x_1 & x_2 & x_3\\
  x_3 & x_1 & x_2 \\
  x_2 & x_3 & x_1\\
\end{pmatrix}
$$
with the determinant
$$
\Theta=(x_1+x_2+x_3)(x_1+\varepsilon x_2 + \varepsilon^2 x_3)(x_1+\varepsilon^2 x_2 + \varepsilon x_3), \quad \varepsilon=e^{\frac{2\pi i}{3}}.
$$

In general, the group determinant $\Theta=\Theta_G \in \mathbb C[x]$ is a homogeneous polynomial of degree $N$ of $$x=(x_g), g \in G=(x_{g_1},\dots, x_{g_N}).$$

Dedekind conjectured that for any finite group $\Theta$ can be factorised over $\mathbb C$ with as many linear factors as the index of the commutator subgroup $G' \subset G$
(in particular, for abelian groups $\Theta$ is always a product of linear factors). 

The nature of the remaining factors was not clear to him, and this is where Frobenius made a crucial contribution by introducing his $k$-characters.

Firstly, he realised that the linear factors of $\Theta$ are simply related to the homomorhisms $\chi: G \to \mathbb C^{\times}$. Namely, he showed that if $\chi$ satisfies
$$
\chi(gh)=\chi(g)\chi(h), \quad g,h \in G
$$
then 
$$
\Psi_{\chi}:=\sum_{g\in G} \chi(g) x_g
$$
is a linear factor of $\Theta$ (see the example above), which explained Dedekind's observation.

Then Frobenius made a crucial observation that the degree $k$ factors must have something to do with the homomorphisms
$$
\varphi: G \rightarrow GL(k, \mathbb C).
$$
This motivated him to develop the theory of the linear representations of finite group, including the foundational results about the characters of irreducible representations. 

By now this is a standard undergraduate material (see e.g. Etingof et al \cite{Etingof} and Serre \cite{Serre}), but it is amazing to see how much of this theory is contained in the pioneering Frobenius paper \cite{Frob1}, including the detailed examples of the symmetry groups of platonic solids and matrix groups over finite fields!

His main result can be stated as follows.
Let $G$ be a finite group and $r$ be the number of conjugacy classes in $G$.

\begin{Theorem} (Frobenius \cite{Frob1, Frob2})
The group determinant $\Theta= \det \, X_G$ has the factorisation of the form
\beq{fact}
\Theta(x)=\prod_{i=1}^r \Psi_i(x)^{\deg \, \Psi_i},
\eeq
where $\Psi_i(x), \, i=1,\dots, r$ are certain irreducible polynomials.
\end{Theorem}

The modern proof is easy (see Buchstaber-Rees \cite{BR-2004} and Etingof et al \cite{Etingof}).
Consider the group algebra $\mathbb C[G]$ and the corresponding regular representation $\rho: G \to End(V), \, V=\mathbb C[G]$ acting by the left multiplication $h \to gh$ on the basis $\{h\}, \, h\in G$ in $V.$

Consider the $End(V)$-valued polynomial
$$
L(x)=\sum_{g \in G}x_g \rho(g).
$$
Its matrix in the basis $\{g_i\}\in V, \, i=1,\dots, N$ is precisely $X_G$ since
$$
L(x)h=\sum_{g \in G}x_g \rho(g)h=\sum_{g \in G}x_g gh=\sum_{g \in G}x_{gh^{-1}} g.
$$
Let $V_i, \, i=1,\dots, r$ be the set of all non-equivalent irreducible representations of $G$, then the regular representation is isomorphic to the direct sum
$$
V=\oplus_{i=1}^{r} (\dim V_i)\, V_i,
$$
which is essentially what Frobenius had shown. In fact, the claim that every irreducible representation appears in the regular representation as often as its dimension was called by Frobenius ``the Fundamental Theorem of the theory of group determinants" and took most of his time to prove (see the Historical interlude 4.11 in \cite{Etingof}).

Now taking the determinant of the corresponding block-diagonal matrix realisation of $L(x)$ we have
$$
\det  L(x)=\prod_{i=1}^r \Psi_i(x)^{\dim \, V_i},
$$
where $\Psi_i(x)=\det L_i(x)$ is the determinant of the corresponding block \cite{Etingof}.

But the most important for us is how Frobenius proposed to compute these polynomials $\Psi_i(x)$ since this is really ingenious and not very well-known.

His main construction was as following. 

Let us call a function $\chi:G \to \mathbb C$ {\it trace-like} if for any $g,h \in G$
$$
\chi(gh)=\chi(hg).
$$
Starting with any trace-like function $\chi: G \to \mathbb C$ Frobenius defines 
the polynomials $\Phi _{k}^{\chi}(x)$
\beq{phin}
\Phi_k^{\chi}=\frac{1}{k! }\sum_{h_1,\dots,h_k\in G}\chi_k(h_1,\dots,h_k)x_{h_1}\dots x_{h_k},
\eeq
where the functions $\chi_k(h_1,\dots, h_k)$ are defined by the recurrence procedure
\beq{recurr}
\chi_{k+1}(h_0, h_1,\ldots, h_k): =\chi(h_0)\chi
_{k}(h_1, h_2,\ldots, h_k )
\eeq
$$
- \sum_{j=1}^k \chi_k(h_1, h_2,\ldots, h_0h_j \ldots, h_k ), \quad h_0, h_1, \dots, h_k \in G.  
$$
with $\chi_1=\chi.$ In particular, $\Phi _{1}^{\chi}=\sum_{h\in G}\chi(h)x_h.$

When $\chi=tr\, \rho(g)$ is the character of some representation $\rho: G \to GL(V)$ then $\chi_k$ is called {\it Frobenius $k$-character} of $V.$

Frobenius wrote the first three of them explicitly: in his notations
\beq{c2}
\chi_2(A,B)=\chi(A)\chi(B)-\chi(AB),
\eeq
$$
\chi_3(A,B,C)=\chi(A)\chi(B)\chi(C)-\chi(A)\chi(BC)-\chi(B)\chi(AC)-\chi(C)\chi(AB)
$$
\beq{c3}
+\chi(ABC)+\chi(ACB)
\eeq
for all $A,B,C \in G$ (see formula (14) in section 3 of \cite{Frob2}).

Since the initial function $\chi$ is trace-like, the functions $\chi_k:G^k \to \mathbb C$ are symmetric functions of group elements. This follows from the following formula for $\chi_n$ presented by Frobenius as a sum over symmetric group $S_k.$

Each permutation $\sigma \in S_k$ can be uniquely represented as a product of disjoint cycles:
$\sigma =\gamma _{1}\dots \gamma
_{s}.$ For a cycle $\gamma =(i_{1} i_2 \ldots i_{m})$ define $$\chi_{\gamma
}(h_{1},\ldots , h_{k})=\chi (h_{i_{1}}h_{i_{2}}\ldots h_{i_{m}}),$$
which is well-defined for every cycle because of the trace-like property of $\chi.$

Frobenius showed that $\chi_k$ can be equivalently written as
\beq{cyclic}
 \chi _{k}(h_{1},h_2,\ldots, h_{k})=
\sum\limits_{\sigma \in S_{k}} \varepsilon(\sigma)\chi_{\gamma
_{1}}(h_{1},\ldots, h_{k}) \ldots \chi_{\gamma _{s}}(h_{1},\ldots, h_{k}),
\eeq
 where
$\varepsilon(\sigma)$ is the sign of the permutation $\sigma$ (see formula (17) in section 3 of \cite{Frob2}).
For $k=3$ one can clearly see this in formula (\ref{c3}).

A remarkable discovery of Frobenius is the following

\begin{Theorem} (Frobenius \cite{Frob2})
Let $V$ be an irreducible representation of $G$ of dimension $n$ and $\chi$ be its character.

Then the corresponding $\chi_j \equiv 0$ for all $j>n$ and $\Phi_n^\chi=\Psi_V$ is the corresponding irreducible factor of the group determinant in (\ref{fact}).
\end{Theorem}

%The functions $\chi$ such that the corresponding $\chi_{k+1}\equiv 0$ (and hence $\chi_j \equiv 0$ for all $j>k$) but $\chi_k\neq 0,$ Frobenius called the characters of group $G$ of $k$-th degree (nowadays usually called {\it Frobenius $k$-characters}).

Frobenius showed that in that case  
\beq{ce}
\chi_{j+1}(e,h_1,\dots,h_j)=(k-j)\chi_j(h_1,\dots, h_j), \quad \chi(e)=\chi_1(e)=k
\eeq
and the corresponding polynomials $\Phi_j(x)=\Phi_j^{\chi}(x)$ appear as the coefficients in the following expansion
\beq{phi}
\Phi_k(x+u\epsilon)=u^k+\Phi_1(x)u^{k-1}+\Phi_2(x) u^{k-2}+\dots + \Phi_k(x),
\eeq
where in Frobenius's notations $$x+u\epsilon=(x_e+u, x_{g_2},\dots, x_{g_N})$$ (see formula (2) in Section 3 of \cite{Frob2}).

Note that the condition that $\chi_{k+1}\equiv 0$ holds for the character $\chi$ of every $k$-dimensional representation of $G$, not necessarily irreducible, and for every group $G$, not necessarily finite.

Formula (\ref{phi}) helps us to reveal the nature of the $k$-characters. 

Let $\chi_1=\chi$ be the character of a $n$-dimensional irreducible representation $\rho: G \to GL(V)$ and consider the function (cf. Atiyah \cite{Atiyah})
\beq{lambda}
\Lambda_u(g)=\det (uI + \rho(g)), \,\, g \in G
\eeq
depending on a parameter $u.$

Following Frobenius \cite{Frob2} consider the restrictions of Frobenius characters on the diagonal:
\beq{diag}
\theta_k^\rho(g)=\frac{1}{k!}\chi_k(g,\dots,g), \,\, g \in G.
\eeq

\begin{prop}
The function $\Lambda_u(g)$ is the generating function of the diagonal Frobenius characters:
\beq{phi2}
\det (uI + \rho(g))=u^n+\theta_1^\rho(g)u^{n-1}+\theta_2^\rho(g) u^{n-2}+\dots + \theta_n^\rho(g).
\eeq
This means that the diagonal Frobenius $k$-characters of representation $\rho$ are simply the usual characters of $k$-th exterior power of $\rho$:
\beq{ext}
\theta_k^\rho=\chi(\Lambda^k \rho).
\eeq
\end{prop}

The proof can be essentially extracted from Frobenius, who used the Waring relation between the elementary symmetric polynomials $e_1,\dots, e_n$ and Newton power sums $s_k=\sum x_i^k$, which can be written in the form
\beq{en}
n! e_n= det   \left(
\begin{array}{cccccc}
s_1 & 1 & 0 & 0& \dots &0\\
s_2 & s_1 & 2  & 0& \dots  &0 \\
\vdots &\vdots & \vdots& \ddots& & \vdots \\
\vdots &&&&&\vdots\\
s_{n-1} &s_{n-2} & s_{n-3} & \dots & s_1
& {n-1 } \\
s_n &s_{n-1} & s_{n-2}  &  \dots &s_2
&s_1\\
\end{array} \right).
\eeq

The observation that the diagonal Frobenius $k$-character is the character of $k$-th exterior power might be new, but it does not characterise the $k$-characters since one cannot reconstruct $\chi$ from its diagonal version $\bar\chi.$ 

However, as it was pointed out by Buchstaber and Rees, the diagonal version of the corresponding $n$-homomorphism $\Phi_n: \mathbb CG \to \mathbb C$ uniquely determines $\Phi_n$ by polarisation, see \cite{BR-2002, BR-2004}.

It is interesting that Atiyah used similar ideas, and in particular function $\Lambda_u(g)$, to describe the filtration in the Grothendieck ring $R(G)$ of the representations of $G$, inspired by the construction of the Chern classes in complex $K$-theory for complex vector bundles, see section 12 in \cite{Atiyah}.

Let us discuss now what happens for Frobenius $k$-characters in the reducible case.

Let $\rho=\rho_1\oplus \rho_2$  be the sum of two representations of $G$, then a simple analysis of Frobenius recursion leads to the following result.

\begin{prop}
Frobenius $k$-character $\chi_k^\rho$ of the reducible representation 
$\rho=\rho_1\oplus \rho_2$ can be expressed via Frobenius characters of $\rho_1$ and $\rho_2$ by the formula
$$
\chi_k^\rho(h_1,\ldots,h_k) = \chi_k^{\rho_1}(h_1,\ldots,h_k) + \chi_k^{\rho_2}(h_1,\ldots,h_k)
$$
\beq{sum}
+\sum_{m=1}^{k-1}\sum_\sigma \chi_{m}^{\rho_1}(h_{\sigma(1)},\ldots, h_{\sigma(m)})
\chi_{k-m}^{\rho_2}(h_{\sigma(m+1)},\ldots, h_{\sigma(k)}),
\eeq
where $\sigma$ runs through all $(m,k-m)$-shuffles, which are permutations from $S_k$ such that $\sigma(1)< \ldots <\sigma(m)$ and $\sigma(m+1)< \ldots <\sigma(k).$ 

The same formula holds for $\chi=\chi^{(1)} + \chi^{(2)}$
for any trace-like functions $\chi^{(1)}, \chi^{(2)}: G \to \mathbb C$. \end{prop}

The number of $(m,k-m)$-shuffles is $k \choose m$, so the sum in the right hand side has $2^k$ terms.

In particular, we have 
\begin{align*}
\chi_1^\rho(h) &= \chi_1^{\rho_1}(h) + \chi_1^{\rho_2}(h), \\
\chi^\rho_2(h_1,h_2) &= \chi^{\rho_1}_2(h_1,h_2)+\chi_1^{\rho_1}(h_1)\chi_1^{\rho_2}(h_2)+\chi_1^{\rho_1}(h_2)\chi_1^{\rho_2}(h_1)+\chi^{\rho_2}_2(h_1,h_2),\\
\chi^\rho_3(h_1,h_2,h_3) &= \chi^{\rho_1}_3(h_1,h_2,h_3)+\chi^{\rho_1}_2(h_1,h_2)\chi_1^{\rho_2}(h_3)+\chi^{\rho_1}_2(h_1,h_3)\chi_1^{\rho_2}(h_2)\\
&+\chi^{\rho_1}_2(h_2,h_3)\chi_1^{\rho_2}(h_1) + \chi_1^{\rho_1}(h_1)\chi^{\rho_2}_2(h_2,h_3)+\chi_1^{\rho_1}(h_2)\chi^{\rho_2}_2(h_1,h_3)\\
&+\chi_1^{\rho_1}(h_3)\chi^{\rho_2}_2(h_1,h_2)+\chi^{\rho_2}_3(h_1,h_2,h_3).
\end{align*}

%$(k-m,m)$-{\it перетаcовки} (shuffle - перетасовка).

For the diagonal versions (\ref{diag})
%$\theta_k^\rho : G \to \mathbb{C}$
%$$
% \theta_k^\rho(h)=\frac{1}{k!}\chi_k^\rho(h,\ldots,h)
%$$
we have the Whitney formula, which is well-known in topology:
%\begin{align*}
%\widehat{\chi}_1^\rho &= \widehat{\chi}_1^{\rho_1} + \widehat{\chi}_1^{\rho_2}; \\
%\widehat{\chi}_2^\rho &= \widehat{\chi}_2^{\rho_1} + \widehat{\chi}_1^{\rho_1} \widehat{\chi}_1^{\rho_2} + \widehat{\chi}_2^{\rho_2}; \\
%\widehat{\chi}_3^\rho &= \widehat{\chi}_3^{\rho_1} + \widehat{\chi}_2^{\rho_1} \widehat{\chi}_1^{\rho_2} + \widehat{\chi}_1^{\rho_1} \widehat{\chi}_2^{\rho_2} + \widehat{\chi}_3^{\rho_2}.
%\end{align*}
%General formula
\beq{whit}
\theta_k^\rho(g) = \sum_{m=0}^k \theta_{m}^{\rho_1}(g) \theta_{k-m}^{\rho_2}(g), \quad g\in G,
\eeq
where $\theta_0^\rho(g) \equiv  1.$

Let us explain now the relation of Frobenius theory to the work of Fricke.

\begin{prop}
Fricke identity (\ref{fricke}) is the corollary of the Frobenius condition $\chi_3(A,B,C)=0$ for the 2-characters.
\end{prop}

Indeed, applying Frobenius theory to canonical 2-dimensional representation of Fuchsian group $G$ we have the relation $\chi_3(A,B,C)=0$ for every $A,B,C\in G$ and $\chi(A)=tr\, A.$ Using the explicit form (\ref{c3}) and substituting into the corresponding relation $A=ab=c,\, B=a^{-1},\, C=b^{-1}$
we have
$$\chi(A)\chi(B)\chi(C)-\chi(A)\chi(BC)-\chi(B)\chi(AC)-\chi(C)\chi(AB)
+\chi(ABC)+\chi(ACB)$$
$$
=\chi(ab)\chi(a^{-1})\chi(b^{-1})-\chi(ab)\chi(a^{-1}b^{-1})-\chi(a^{-1})\chi(a)-\chi(b^{-1})\chi(aba^{-1})
$$
$$+\chi(aba^{-1}b^{-1})+\chi(e)=\chi(a)\chi(b)\chi(c)-\chi(a)^2-\chi(b)^2-\chi(c)^2+\chi(aba^{-1}b^{-1})+2=0,
$$
which coincides with Fricke identity (\ref{fricke}).

We can interpret this calculation as follows: the Fricke identity corresponds to Frobenius relation under the inverse map $A \to A^{-1},$ which establishes an isomorphism between two group structures on $G$:
$m_L: A,B \to AB$ and $m_R: A,B \to BA.$

\section{Modern developments}

The theory of $k$-characters of Frobenius appears to be largely forgotten until 1990s, when Johnson raised the question whether they can be used to characterise the group uniquely.
It turned out that this is indeed the case and it is enough to use only $k$-characters with $k\leq 3$ (see \cite{HJ, Joh}). Formanek and Sibley  proved that the group determinant also uniquely determines the finite group \cite{Formanek}.

Buchstaber and Rees \cite{BR-2004, BR-2008} developed a closely related theory of the $n$-homomorphisms and $n$-Hopf algebras, motivated by the theory of multi-valued groups and of symmetric products of topological spaces.
This allows them to give a short proof of Hoenke-Johnson theorem about characterisation of the group by Frobenius $k$-characters with $k\leq 3$ (see Sectrion 4 in \cite{BR-2004}).
They have used the notion of Frobenius algebra, which also plays the key role in Dubrovin's theory of Frobenius manifolds.

The works of Markov and Fricke were much more fortunate. First of all, it attracted the attention of Frobenius, who in 1913 wrote a paper "On Markov numbers" \cite{Frob3}, where he studied them in more detail and conjectured that the Markov triple $(x,y,z)$ is uniquely determined by its largest element (the celebrated {\it Uniqueness Conjecture}, which is still open, see Aigner \cite{Aigner}).

As we have already mentioned, in 1950s Gorshkov and Cohn linked Markov triples with the length of simple geodesics on  punctured torus with special hyperbolic metric \cite{Cohn, Gorshkov}. The Fricke work \cite{Fricke}  is now well recognised as one of the foundational for the future Teichm\"uller theory.

Then it came the turn of algebraic geometry, where Markov triples surprisingly appeared in a number of different connections. 

In 1988 Rudakov \cite{Rud} showed the Markov triples are precisely the ranks of the vector bundles in the exceptional collections on projective plane $\mathbb P^2.$ 
More recently, Hacking and Prokhorov \cite{HP} showed that the smoothable del Pezzo surfaces with quotient singularities are the weigthed projective planes with weights given by the squares $(x^2,y^2,z^2)$ of the Markov numbers.

But the most remarkable link from our point of view was discovered by Dubrovin in 1994 in relation with his theory of Frobenius manifolds, enumerative geometry and quantum cohomology of $\mathbb P^2.$

Let $N_d$ be the numbers of rational curves in $\mathbb P^2$ passing through generic $3d-1$ points. Kontsevich and Manin \cite{KM} showed that certain generating function of them satisfies the associativity equation and thus the Painlev\'e-VI equation with very special parameters. Dubrovin \cite{D1} showed that the initial data for the corresponding solution is nothing other but the Markov orbit! 

Finally we would like to mention very interesting recent work by Bourgain, Gamburd and Sarnak \cite{BGS}, who argued in favour of the {\it Strong Approximation Conjecture} claiming that all the non-zero solutions of the Markov equation modulo prime $p$ are coming from the Markov triples.

It is interesting that the solutions modulo $p$ of a version of Markov equation
$$
x^2+y^2+z^2-2xyz=1
$$
appears already in Frobenius's study of the characters of the finite group $SL_2(\mathbb F_p)$ (see section 10 in \cite{Frob2}).
This equation was studied later by Mordell \cite{Mordell} and used by Zagier \cite{Zagier} to study the asymptotics of Markov numbers.

\section{Acknowledgements}

We are very grateful to Leonid Chekhov, Vladimir Fock and Leon Takhtajan for very stimulating discussions of the role of the Markov triples in the theory of Teichm\"uller spaces, and to John McKay from whom we first learnt about work of Frobenius on group determinants. 

We are particularly indebted to our dear friend, the late Boris Dubrovin, who inspired our interest in Markov triples by explaining their remarkable relation to enumerative geometry and Painlev\`e equations.

The work of the first author was partially supported by the Russian Ministry of Education and Science RF (project N1.13560.2019/13.1).

\end{document}